\documentclass[12pt,french]{amsart}
\usepackage[french]{babel}
\usepackage{amsmath,amssymb,latexsym}
\usepackage[latin1]{inputenc}

\def\im{\mathrm{im}}

\def\ncl#1{\overline#1^{\scriptscriptstyle n}}

\theoremstyle{plain}
\newtheorem{theorem}{Th\'eor\`eme}[section]
\newtheorem{proposition}[theorem]{Proposition}
\newtheorem{fact}[theorem]{Fait}
\newtheorem{lemma}[theorem]{Lemme}
\newtheorem{corollary}[theorem]{Corollaire}

\theoremstyle{definition}
\newtheorem{definition}[theorem]{D\'efinition}
\newtheorem{remark}[theorem]{Remarque}
\newtheorem{expl}[theorem]{Exemple}
\newtheorem{frage}[theorem]{Question}

\def\bsp{\begin{expl}}
\def\ebsp{\end{expl}}
\def\beh{\begin{claim}}
\def\ebeh{\end{claim}}
\def\defn{\begin{definition}}
\def\edefn{\end{definition}}
\def\satz{\begin{theorem}}
\def\esatz{\end{theorem}}
\def\tats{\begin{fact}}
\def\etats{\end{fact}}
\def\kor{\begin{corollary}}
\def\ekor{\end{corollary}}
\def\lmm{\begin{lemma}}
\def\elmm{\end{lemma}}
\def\bem{\begin{remark}}
\def\ebem{\end{remark}}
\def\bew{\begin{proof}}
\def\ebew{\end{proof}}
\def\satzli{\begin{proposition}}
\def\esatzli{\end{proposition}}
\def\frg{\begin{frage}}
\def\efrg{\end{frage}}

\begin{document}
\title[Commentaires sur un résultat
d'Olivier Frecon]{Commentaires sur un résultat\\
d'Olivier Frecon}
\author{Bruno Poizat \& Frank O. Wagner}
\address{Universit\'e de Lyon; CNRS; Universit\'e Claude Bernard Lyon 1; 
Institut Camille Jordan UMR5208, 43 bd du 11 novembre 1918, 69622 Villeurbanne 
Cedex, France}
\email{poizat@math.univ-lyon1.fr}
\email{wagner@math.univ-lyon1.fr}

\date{\today}
\thanks{Partially supported by ANR-13-BS01-0006 ValCoMo}
\begin{abstract}Il n'y a pas de groupe malsain de rang $2n+1$ avec un borel 
abélien de rang $n$. En particulier, le théorème de Frécon en découle~: il n'y 
a pas de mauvais groupe de rang $3$.\end{abstract}
\maketitle

\section{Introduction} Dans un preprint récent, Olivier Frécon a donné une 
réponse négative à la question de l'existence d'un mauvais groupe de rang de 
Morley $3$. Sa démonstration fabrique de manière combinatoire et inattendue un 
ensemble $X$ (qualifié ici {\em de Frécon}) qui fournit un automorphisme 
involutif du mauvais groupe, ce qui mène à une contradiction. Nous analysons 
ici sa démonstration et la généralisons aux mauvais groupes de rang $2n+1$ avec 
un borel abélien de rang $n$.

\noindent{\em Notation.} Comme d'habitude, on note $g^h=h^{-1}gh$ et 
$[g,h]=g^{-1}h^{-1}gh$.
\section{Groupes malsains et mauvais groupes}
\defn Un groupe {\em  malsain} est un groupe connexe $G$ de rang de Morley fini 
avec un sous-groupe définissable propre $B$, dont on appellera {\em borels} les 
conjugués, vérifiant les deux choses suivantes~:\begin{enumerate}
\item $B$ est {\em malnormal}, c'est-à-dire que $B\cap B^a = \{1\}$ si $a\notin 
B$~;
\item $G$ est la réunion de ses borels.\end{enumerate}\edefn
En particulier, $B$ est autonormalisant et genereux dans $G$~; comme $G$ n'a 
pas d'involution (lemme \ref{l3}), son existence contredirait la Conjecture 
d'Algébricité. On aimerait bien montrer qu'il n'existe pas de groupe malsain.
\lmm\label{l1} Soit $G$ un groupe malsain. Alors ses borels sont infinis et 
connexes, et tout sous-groupe fini de 
G est contenu dans un borel. De plus, $G$ est sans involution.\elmm
\bew Considérons $a\not= 1$ dans un borel $B$. Si $g\in G$ et $a^g\in B$, alors 
$a^g\in B\cap B^g$, d'où $B=B^g$ et $g\in N_G(B)=B$. En particulier $C_G(a)\le 
B$. Si les borels sont finis, tout élément non-trivial a un centralisateur fini 
et une classe de conjugaison générique dans $G$~; comme $G$ est connexe, tous 
les éléments non-triviaux de $G$ sont conjugués, ce qui est impossible pour un 
groupe stable.

Les borels sont donc infinis~; s'ils n'étaient pas connexes, on obtiendrait 
deux ensembles génériques disjoints: la réunion des conjugués de ~$B^0$, et la 
réunion des conjugués de $B\setminus B^0$.

Soit $F$ un sous-groupe fini non-trivial de $G$., et $H_1$,\ldots,$H_n$ les 
intersections non-triviales de $F$ avec les borels de $G$. Alors les $H_i$ sont 
autonormalisants dans $F$. Si $ H_1$ est strictement inclus dans $G$, le nombre 
$n$ de points de la réunion des $F$-conjugués de $H_1$ est 
$$n=|F:H_1|\cdot(|H_1| -1) + 1 = |F| - |F:H_1| + 1.$$
Comme $H_1$ est d'indice au moins deux dans $F$, on a $n<|F|$ et il y a un 
$H_i$ non-conjugué de $H_1$.
D'autre côté, comme $H_1$ est d'indice au plus $|F|/2$, on a $n\ge1+|F|/2$ et 
il 
n'y a pas d'assez de place pour les conjugués de $H_i$. Ainsi $H_1=F$, et $F$ 
est contenu dans un seul borel.

Enfin, si $i$ et $j$ sont des involutions prises dans des borels différents, 
elles inversent toutes deux $ij$, et normalisent son borel, c'est-à-dire lui 
appartiennent, ce qui est absurde.
\ebew
\lmm\label{l0} Soit $G$ un groupe malsain. Alors tout sous-groupe connexe 
définissable $H$ qui n'est pas inclus dans un borel est encore malsain, et ses 
borels sont les traces des borels de $G$ dans $H$.\elmm
\bew Soient $B_0$,\ldots,$B_i$,\ldots les traces non-triviales dans $H$ des 
borels de $G$. Notons d'abord que chaque $B_i$ est autonormalisant dans $H$. 
Comme $B_0$ est en plus disjoint de ses conjugués dans $H$ en dehors de $1$, la 
réunion
de ces conjugués est une partie générique de $H$. Comme $H$ est connexe, tous 
les $B_i$ sont conjugués dans $H$, et $H$ est malsain.\ebew
\lmm Soit $G$ un groupe malsain avec borel $B$, et $N$ un sous-groupe 
définissable normal non-trivial de $G$. Alors $G=NB$.\elmm
\bew Si $N$ était fini, il serait central et normalisait tous les borels, ce 
qui est absurde. Donc $N$ est infini, et on peut supposer $N$ connexe et 
$H=N\cap B$ non-trivial. Alors $N$ est un groupe malsain avec borel $H$. Soit 
$g\in G$. Alors $H^g$ est un borel de $N$, et il y a $n\in N$ avec $H^n=H^g$, 
d'où $gn^{-1}\in N_G(B)=B$.\ebew
\lmm\label{l2.5} Soit $G$ un groupe malsain avec borel $B$, et $H$ un 
sous-groupe définissable de $G$ contenant $B$. Alors $H$ est connexe, et $G$ 
est aussi malsain de borel $H$.\elmm
\bew Supposons d'abord $H$ connexe. Comme $H$ contient $B$ et les conjugués de 
$B$ recouvrent $G$, les conjugués de $H$ recouvrent $G$.

Soit $H^g$ un conjugué de $H$ avec intersection $H\cap H^g$ non-triviale. Soit 
$1\not=x\in H\cap H^g$~; comme $H$ est malsain de borel $B$, on peut supposer 
$x\in B$. Alors $B$ est contenu dans $H$ et dans $H^g$. Ainsi $B$ et $B^g$ sont 
deux borels de $H^g$, et il y a $h\in H^g$ avec $B^h=B^g$, d'où $gh^{-1}\in 
N_G(B)=B\le H^g$. Ainsi $g\in H^g$ et $H=H^g$, et $H$ est malnormal.

Si $H$ n'est pas connexe, alors $H^0$ contient encore $B$ par connexité de ce 
dernier. Alors $G$ est malsain de borel $H^0$. En particulier $H^0$ est 
malnormal, et $H\le N_G(H^0)=H^0$. Ainsi $H$ est connexe.\ebew
\lmm\label{l2.3} Soit $G$ un groupe de rang de Morley fini sans involutions. 
Alors chaque élément $g$ de $G$ a une unique racine carrée, qui est contenu 
dans 
tout sous-groupe définissable contenant $g$. De plus, un élémént non-trivial 
n'est pas conjugué à son inverse.\elmm
\bew Soient $g\not=1$ un élément de $G$, et $A$ le plus petit sous-groupe 
définis\-sable le
contenant. Alors $A$ est commutatif, et comme il n'a pas d'involutions,
l'homomorphisme de $A$ dans $A$ qui à $x$ associe $x^2$ est injectif, et donc
bijectif puisqu'il conserve rang et degré de Morley. Ainsi tout élément $a\in 
A$ possède une
unique racine carrée $u\in A$~; si $v$ est une autre racine carrée de $a$, elle 
commute avec $a$, normalise $A$, et commute avec l'unique racine carrée de $a$ 
dans $A$, soit $u$. Donc $(uv^{-1})^2 = aa^{-1} = 1$ et $u = v$ puisqu'il n'y a 
pas d'involutions.

Enfin, si $g^h=g^{-1}$, alors $h^2\in C_G(g)$, et donc $h\in C_G(g)$, puisque 
dans un groupe de rang de Morley fini les involutions se relèvent. Ainsi 
$g=g^{-1}$ et $g=1$.\ebew

\lmm\label{l2.4} Soit $G$ un groupe de rang de Morley fini sans involutions, et 
$\sigma$ un automorphisme involutif de $G$. Soit $F=\{g\in 
G:\sigma(g)=g\}$ et $I=\{g\in G:\sigma(g)=g^{-1}\}$. Alors $G=F\cdot 
I=I\cdot F$, avec décomposition unique~; de plus, un point non-trivial de 
$F$ n'a pas de conjugué dans $I$.\elmm
\bew Soit $g\in G$ et considérons $h= g^{-1}\sigma(g)$. Alors 
$\sigma(h)=h^{-1}$ et $h\in I$. Soit $\ell$ la racine carrée unique de $h$. 
Alors $\ell^{-1}$ est la racine carrée unique de $h^{-1}=\sigma(h)$, d'où 
$\sigma(\ell)=\ell^{-1}$ et $\ell\in I$. Ainsi
$$\sigma(g\ell)=\sigma(g)\sigma(\ell)=gg^{-1}\sigma(g)\ell^{-1}
=gh\ell^{-1}=g\ell^2\ell^ { -1 }=g\ell$$
et $g\ell\in F$. Donc $g\in F\cdot I$ et $G=F\cdot I$. En prenant des inverses, 
on obtient $G=I\cdot F$.

Si $g=if=i'f'$ avec $i,i'\in I$ et $f,f'\in F$, alors $i^{-1}i'=ff'^{-1}\in F$, 
et
$$i^{-1}i'=\sigma(i^{-1}i')=\sigma(i^{-1})\sigma(i)=ii'^{-1},$$
d'où $i'^2=i^2$. Alors $i=i'$ et la décomposition est unique.

Enfin, soit $f\in F$ non-trivial, et $g\in G$ tel que $f^g\in I$. Alors $g=f'i$ 
pour un $i\in I$ et $f'\in F$, et $f^{f'}\in F$ est non-trivial. On peut donc 
supposer $g=i\in I$. Alors 
$$i^{-1}fi=f^i=\sigma(f^i)=\sigma(i^{-1}fi)=ifi^{-1}$$
et $fi^2=i^2f$. Donc $i^2\in C_G(f)$; d'après le lemme \ref{l2.3} on a aussi 
$i\in C_G(f)$, et $f^i=f\in F\cap I=\{1\}$, une contradiction. 
\ebew
\bem Dans le lemme \ref{l2.4} on ne suppose pas que l'automorphisme involutif 
$\sigma$ soit définissable.\ebem

La famille des borels d'un groupe malsain n'est pas a priori uniquement 
déterminée.
On appelle un groupe malsain {\em canonique} si ses borels sont ses 
sous-groupes définissables connexes maximaux.

\lmm Soit $G$ un groupe malsain. Un sous-groupe $H$ définissable connexe de 
rang minimal qui n'est pas inclus dans un borel est malsain canonique.\elmm
\bew $H$ est malsain d'après le lemme \ref{l0}. Soit $F$ un sous-groupe propre 
définissable de $H$. Si $F$ est fini, il est
dans un borel d'après le Lemme \ref{l1}. S'il est infini, $F^0$ est dans un 
borel $B$ par minimalité de $H$, et $F$ normalise $F^0$, donc $B$, et $F\le 
B$.\ebew

\satzli\label{p2.10} Un groupe malsain canonique $G$ n'a pas d'automorphisme 
involutif non-trivial.\esatzli
\bew Soit $\sigma$ un automorphisme involutif de $G$. On pose 
$F=\mbox{Fix}_G(\sigma)$ et $I=\{g\in G:\sigma(g)=g^{-1}\}$. Alors $F^0$ est 
contenu dans un borel $B$~; puisque $B$ est malnormal, $F\le B$. Comme 
$G=F\cdot I$, il y a un autre borel $B'$ qui contient un élément $i$ non-trivial 
de $I$. Alors $i^{-1}\in B'\cap\sigma(B')$ et $\sigma(B')=B'$. Ainsi $\sigma$ 
est un automorphisme involutif de $B'$, et $B'=(I\cap B')\cdot(F\cap B')$ 
d'après le lemme \ref{l2.4}. Or, $B\cap F=\{1\}$, et $B'=I\cap B'$ est inversé 
par $\sigma$. Mais $B$ se conjugue dans $B'$, ce qui contredit le lemme 
\ref{l2.4}.\ebew

Rappelons qu'un mauvais groupe est un groupe connexe non-nilpotent de rang de 
Morley fini dont tous les sous-groupes définissables connexes résolubles sont 
nilpotents. Le premier mauvais groupe est apparu dans un ancien article de 
Cherlin \cite{Ch}, comme un groupe simple de rang de Morley trois qui n'est pas 
de la forme PSL$_2(K)$; c'est Nesin \cite{Ne} qui a montré leur malignité au 
sens ci-dessus, c'est-à-dire l'absence d'involutions, en s'appuyant sur un 
théorème de Bachman. Il a été ensuite remarqué par Corredor \cite{Co} et 
Borovik-Poizat \cite{BP} que les mauvais groupes produisaient des groupes 
malsains canoniques. D'autres situations semblant paradoxales, comme un groupe 
minimal avec un groupe
d'automorphismes définissable non-commutatif, ou bien un groupe 
définis\-sablement linéaire contredisant la Conjecture d'Algébricité, en 
produisent également. Frécon \cite{Fr} a réussi dernièrement à montrer que les 
mauvais groupes de rang
de Morley trois n'existent pas~; les deux sections suivantes essaient de 
généraliser ses arguments.

\section{Ensembles de Frécon}
Soit $G$ un groupe malsain, avec un borel $B$. Si $u\not=1$ est dans $G$, on 
notera $B(u)$ l'unique borel contenant $u$.
\defn On appelle {\em droite} un translaté bilatère de $B$, c'est-à-dire un 
ensemble de la forme $uBv = uv\cdot B^v$, qui est donc un translaté à gauche du 
borel $B^v$, qu'on appelle sa {\em  direction}. L'ensemble des droites est 
not\'e $\Lambda$.\edefn
Pour une partie $X$ de $G$ on pose $X^{-1}=\{x^{-1}:x\in X\}$. Alors l'ensemble 
des droites est préservé par translation à gauche et à droite, et aussi par 
inversion.
Deux droites distinctes de même direction sont dites {\em parallèles}~; deux 
droites parallèles sont disjointes. Notons que le passage des translatés à 
gauche aux translatés à droite ne conserve pas le parallélisme.

Si deux droites distinctes se coupent, elles forment un translaté modulo
l'intersection de leurs directions respectives, qui est réduite à l'élément 
neutre~;
par conséquent, deux droites distinctes sont disjointes ou se coupent en un seul
point. Si $u$ et $v$ sont deux points distincts, ils sont contenus tous les 
deux dans la droite $uB(u^{-1}v)$. C'est donc la seule droite qui passe par ces 
deux points~; on la note $\ell(u,v)$.

Nous dirons qu'une partie définissable de G est {\em convexe} si avec deux
points elle contient toute la droite qui les joint. Un ensemble convexe fini 
est un
point~; une droite, ou plus généralement un translaté d'un sous-groupe connexe 
contenant un borel est convexe. En fait, nous allons voir que ce sont les seuls 
ensembles convexes. Ce sera conséquence d'un théorème plus fort, concernant des 
ensembles plus généraux que les
ensembles convexes, et qui --- nous l'espérons --- seront plus faciles à faire 
paraître
dans les démonstrations d'inexistence.

Nous dirons qu'une droite $\ell$ est {\em génériquement incluse} dans l'ensemble
définissable $X$ si $RM(X\cap\ell)=RM(\ell)$~; comme la dimension est 
définissable,
les droites génériquement incluses dans $X$ forment une famille définissable 
$\Lambda(X)$.
Notons que si $RM(X)<RM(B)$ alors $\Lambda(X)=\emptyset$, et si $RM(X\triangle 
Y)<RM(B)$, alors $\Lambda(X)=\Lambda(Y)$.

Un ensemble définissable $X$ est {\em quasi-convexe} si pour tout $x,y\in X$ 
distincts, $\ell(x,y)$ est génériquement incluse dans $X$.
Un ensemble convexe est bien sur quasi-convexe.
\lmm Un ensemble $X$ quasi-convexe est de degré de Morley $1$.\elmm
\bew On considère $$Z=\{(x,y,\ell):x,y\in X\cap\ell,\ 
\ell\in\Lambda(X)\}\subseteq X^2\times\Lambda(X).$$
Alors $(x,y,\ell)$ est générique dans $Z$ si et seulement si $x,y$ sont 
génériques indépendants dans $X$ et $\ell=\ell(x,y)$, ou encore si $\ell$ est 
générique dans $\Lambda(X)$ et $x,y$ sont génériques indépendants dans $\ell$. 
Mais $\ell$ n'a qu'un seul type générique, ce qui implique que $x$ et $y$ on 
même type (fort) sur $\emptyset$, et $DM(X)=1$.\ebew
\lmm\label{l3} $RM(\Lambda(X))\le 2\,RM(X)-2\,RM(B)$.\elmm
\bew On considère l'application $f:X^2\to\Lambda$ donnée par 
$f(x,y)=\ell(x,y)$. Si $\ell\in\Lambda(X)$, alors le fibre sur $\ell$ est de 
rang $2\,RM(B)$~; le lemme en découle par additivité du rang de Morley.\ebew
\bem Pour $RM(X)\ge RM(B)$ le maximum est atteint\begin{itemize}
\item si $X$ est quasi-convexe infini~;
\item pour $RM(X)=RM(B)$ si et seulement si $X$ contient générique\-ment une 
droite~;
\item si $X$ est générique dans $G$.\end{itemize}
Plus généralement, le maximum est atteint si $X$ contient génériquement un 
translaté $C$ d'un sous-groupe définissable contenant un borel, puisque tous 
les points géné\-riques d'une droite générique sont génériques dans $C$.
Nous allons voir (théorème \ref{t1}) que ce sont les seuls cas possibles.\ebem

Pour une partie $X\subseteq G$ et un élémént $x\in G$ on notera 
$$\Lambda_x(X)=\{\ell\in\Lambda(X):x\in\ell\}$$
l'ensemble des droites presque contenues dans $X$ qui passent par $x$.
\defn Un point $x\in G$ est {\em $n$-interne} à $X$ si $RM(\Lambda_x(X))\ge 
n-RM(B)$~; l'ensemble des points $n$-internes à $X$ est la {\em $n$-clôture} de 
$X$, notée $\ncl X$.\edefn
Comme les droites de $\Lambda_x(X)$ ne s'intersectent qu'en $x$ et sont 
géné\-riquement contenues dans $X$, on a pour $x\in\ncl X$ que
$$RM(X)\ge RM(\bigcup\Lambda_x(X))=RM(\Lambda_x(X))+RM(B)\ge n.$$
En particulier, si $RM(X)<n$ alors $\ncl X=\emptyset$, et si $\ncl 
X\not=\emptyset$ pour un entier $n$, alors $RM(X)\ge RM(B)$.
\lmm\label{l4} Si $RM(X\triangle Y)<n$, alors $\ncl X=\ncl Y$.\elmm
\bew Soit $x\in\ncl X\setminus\ncl Y$. Alors 
$$\bigcup_{\ell\in\Lambda_x(X)}(\ell\cap X)\setminus Y\subseteq X\setminus Y$$
est de rang de Morley au moins $n-RM(B)+RM(B)=n$, une contradiction.\ebew
\lmm\label{l5} Pour une partie $X$ de $G$ et $RM(B)\le n\le RM(X)$ on a
$$RM(\ncl X)+n\le RM(\Lambda(X))+2\,RM(B)\le\max_n\{RM(X\cap\ncl X)+n\}.$$\elmm
\bew Considérons $Z=\{(x,\ell)\in G\times\Lambda(X):x\in\ell\}$ et 
$Z'=\{(x,\ell)\in Z:x\in X\}$. Alors
$$\begin{aligned}RM(Z)&=RM(\Lambda(X))+RM(B)=RM(Z'),\qquad\mbox{et}\\
RM(Z)&=\max_{RM(B)\le n\le RM(X)}\{RM(\ncl X\setminus\overline 
X^{\scriptscriptstyle{n+1}})+n-RM(B)\},\\
RM(Z')&=\max_{RM(B)\le n\le RM(X)}\{RM((\ncl X\setminus\overline 
X^{\scriptscriptstyle{n+1}})\cap X)+n-RM(B)\}.\end{aligned}$$
Le lemme en découle.\ebew

\lmm Soit $X$ une partie de $G$ de rang de Morley $n$. Alors 
$RM(\Lambda(X))=2n-2\,RM(B)$ si et seulement si $RM(\ncl X)=n$. Dans ce cas, si 
$X$ est de degré de Morley $1$, alors $\Lambda(X)$ et $\ncl X$ le sont aussi, 
et $RM(X\triangle\ncl X)<n$.\elmm
\bew Si $RM(\Lambda(X))=2n-2\,RM(B)$, l'ensemble
$$\{(x,y)\in X^2:\ell(x,y)\in\Lambda(X)\}$$
est générique dans $X^2$. Ainsi pour $x,y$ génériques indépendants de $X$ on a 
$y\in\Lambda_x(X)$, d'où $RM(\Lambda_x(X))\ge n-RM(B)$ et $x\in\ncl X$. Donc 
$$RM(\ncl X)\ge RM(X\cap\ncl X)\ge n.$$
L'autre direction découle du lemme \ref{l5}.

Si de plus $X$ est de degré de Morley $1$, une droite générique de $\Lambda(X)$ 
est de la forme $\ell(x,y)$ pour $x,y$ génériques indépéndants de $X$. Donc 
$\Lambda(X)$ est de degré de Morley $1$, ainsi que l'ensemble $Z$ de la 
démonstration du lemme \ref{l5}, et aussi $\ncl X$. Mais $RM(X\cap\ncl X)=n$, 
d'où $RM(X\triangle\ncl X)<n$.\ebew
\defn Une partie $X\not=\emptyset$ de $G$ est un {\em $n$-ensemble de Frécon} 
si $X\subseteq\ncl X$. Si en plus $n=RM(X)$ et $X=\ncl X$ est de degré de 
Morley $1$, nous appelons $X$ un {\em ensemble de Frécon}.\edefn
Notons que tout ensemble quasi-convexe infini est un ensemble de Frécon.
\lmm\label{l8} Si $X$ est une partie de $G$ avec $RM(\ncl X)=RM(X)=n$, il y a 
un ensemble de Frécon de rang $n$ génériquement inclus dans $X$.\elmm
\bew  Soit $X$ la réunion disjointe de parties définissables $X_i$ de rang de 
Morley $n$ et degré de Morley $1$. Comme les droites sont de degré de Morley 
$1$, on a $\Lambda(X)=\bigcup_i\Lambda(X_i)$, encore une réunion disjointe. 
Puisque $RM(\ncl X)=n$ on a $RM(\Lambda(X))=2n-2\,RM(B)$, et il y a $i$ tel que 
$RM(\Lambda(X_i))=2n-2\,RM(B)$. Alors $RM(\ncl X_i)=n$ et $DM(\ncl X_i)=1$. De 
plus, $RM(X_i\triangle\ncl X_i)<n$, d'où $\ncl X_i=\ncl{{\ncl X_i}}$, et $\ncl 
X_i$ est un ensemble de Fr\'econ génériquement inclus dans $X$.\ebew

\lmm Soit $X$ un ensemble de Frécon avec $1\in X$. Alors $X=X^{-1}$, et pour 
tout $x\in X$ on a $xX=Xx^{-1}$.\elmm
\bew $X^{-1}$ est encore un ensemble de Frécon. Comme $\Lambda_1(X)$ consiste 
de borels, on a
$\Lambda_1(X)=\Lambda_1(X^{-1})$. Or, $RM(\bigcup\Lambda_1(X))=RM(X)$ et
$$RM(\bigcup\Lambda_1(X)\setminus X)< RM(\Lambda_1(X))+RM(B)= RM(X)~;$$
comme $X$ est de degré de Morley $1$ on a 
$$RM(X\triangle\bigcup\Lambda_1(X))<RM(X)=RM(\bigcup\Lambda_1(X))$$
et d'après le lemme \ref{l4} on a
$$X=\ncl 
X=\ncl{{\bigcup\Lambda_1(X)}}=\ncl{{\bigcup\Lambda_1(X^{-1})}}=\ncl{{X^{-1}}}=X^
{-1}.$$
Pour $x\in X$ le translaté $Xx^{-1}$ est encore un ensemble de Frécon contenant 
$1$. Donc
$$Xx^{-1}=(Xx^{-1})^{-1}=xX^{-1}=xX.\qedhere$$
\ebew

\satz\label{t1} Soit $G$ un groupe malsain. Alors un ensemble de Frécon est un 
translaté d'un sous-groupe définissable contenant un borel. En particulier, il 
est convexe.\esatz
\bew Supposons, pour une contradiction, que $X$ soit un ensemble de Frécon 
infini qui n'est pas un translaté d'un sous-groupe définissable contenant un 
borel. En translatant, on peut supposer que $1\in X$. Soit $G_0$ le plus petit 
sous-groupe définissable contenant $X$. Alors $G_0$ contient génériquement une 
droite $\ell$, et donc le borel $\ell\ell^{-1}$. Donc $G_0$ est connexe d'après 
le lemme \ref{l2.5}, malsain, et $X$ y est propre par hypothèse.

Soit $N=\{g\in G_0:gX=X\}$ le stabilisateur à gauche de $X$~; comme $X=X^{-1}$, 
c'est également le stabilisateur à droite~; on a $N\subseteq X\subset G_0$. On 
considère le sous-groupe $S$ de $G_0^2$ donné par
$$S=\{(x,y)\in G_0^2:xXy^{-1}=X\}.$$
Alors la projection de $S$ sur les deux coordonnées contient $X$, et doit être 
$G_0$ entier. De plus,
$$\{g\in G_0:(g,1)\in S\}=N=\{g\in G_0:(1,g)\in S\}.$$
Comme $N\times\{1\}$ est normal dans $S$, la projection $N$ sur la première 
coordonnée est normale dans la projection de $S$ sur la première coordonnée, 
soit $G_0$. Ainsi $S$ est le graphe d'un automorphisme $\sigma$ de $G_0/N$ tel 
que $\sigma(x)\equiv x^{-1}\mod N$ pour $x\in X$. En particulier $\sigma^2$ 
fixe $X/N$ et donc $G_0/N$, et $\sigma$ est un automorphisme involutif de 
$G_0/N$. Si $G$ est malsain canonique, $G_0=G$ et $N$ est trivial, et on 
conclut avec la proposition \ref{p2.10}. Sinon, 
on notera que l'ensemble des translatés bilatères de $X$ dans $G_0$  est égal à 
l'ensemble des translatés à gauche, ou à droite.

Soit $H=N_{G_0}(X)=\{g\in G_0:X^g=X\}$ le normalisateur de $X$ dans $G_0$, et 
$g\in G_0$. Il y a $h\in G_0$ avec $X^g=Xh$~; soit $u$ la racine carrée de $h$ 
modulo $N$. Comme $1\in X^g$ on a $h\in X^{-1}=X$ et $\sigma$ inverse $h$ 
modulo $N$, et aussi $u$. Alors
$$X^g=Xh=Xu^2=\sigma^{-1}(u)Xu=u^{-1}Xu=X^u,$$
et $gu^{-1}\in H$. Ainsi $G_0=HX'$, où $X'$ est l'ensemble des racines carrées 
des éléments de $X$. Notons que $\sigma(x)\equiv x\mod N$ si et seulement si 
$x\in H$, et $\sigma$ inverse $X'$ modulo $N$.

Soit $g\in G_0$ avec $\sigma(g)\equiv g^{-1}\mod N$, et choisissons $h\in H$ et 
$x\in X'$ avec $hx=g$. Alors
$$x^{-1}h^{-1}=g^{-1}\equiv\sigma(g)=\sigma(hx)=\sigma(h)\sigma(x)\equiv 
hx^{-1}\mod N.$$
Ainsi $h^{-x}\equiv h\mod N$ et $x^2\in C_{G_0}(h/N)$. Comme $G_0$ n'a pas 
d'involutions, on obtient $x\in C_{G_0}(h/N)$ et $h\in N$. Donc $NX'$ est 
l'ensemble des éléments de $G_0$ inversés par $\sigma$ modulo $N$, qui est clos 
par carré et racine carrée. Si $x\in X'$ et $u\in NX'$ avec $u^2=x$, alors 
$u=ny$ avec $n\in N$ et $y\in X'$, et
$$x=u^2=(ny)^2=nn^{y^{-1}}y^2\in NX=X.$$
Ainsi $X'\subseteq X$~; comme $\sigma$ inverse $X$ modulo $N$ on a 
$$NX'\subseteq NX=X\subseteq NX',$$
et $X=NX'$. En particulier $G_0=HX$ implique $N<H$.

Soit $B$ un borel presque contenu dans $X$. Pour $x,y$ génériques 
indé\-pendants de $B$ on a 
$$y^{-1}x^{-1}\equiv\sigma(xy)=\sigma(x)\sigma(y)\equiv x^{-1}y^{-1}\mod N$$
et $B/N$ est abélien, inversé par $\sigma$. Ainsi $B$ est contenu entièrement 
dans $X$~; comme ceci est vrai pour tout translate de $xX$ avec $x\in X$, toute 
droite génériquement contenue dans $X$ y est contenue, et $X$ est convexe.

Soit $h\in H\setminus N$, et $B$ un borel avec $h\in B$, et soit $B^x$ un 
conjugué de $B$ contenu dans $X$. Alors pour $y=x\sigma(x)$ on a
$$By=Bx\sigma(x)=xB^x\sigma(x)\subseteq xX\sigma(x)=X.$$
Ainsi $y\in X$ et $hy\in X$, d'où
$$h\equiv\sigma(h)=\sigma(hy)\sigma(y^{-1})\equiv y^{-1}h^{-1}y=h^{-y}\mod N,$$
et $y^2\in C_{G_0}(h/N)$. Ainsi $y\in C_{G_0}(h/N)$ et $h\in N$, ce qui est 
absurde.\ebew
\kor\label{k1} Soit $G$ un groupe malsain avec un borel $B$. Si $X$ est une 
partie de $G$ de rang $n>0$ avec $RM(\ncl X)=n$, alors $X$ contient 
génériquement un translaté d'un sous-groupe connexe de rang $n$ qui contient un 
borel.\ekor
\bew Il y a un ensemble de Frécon de rang $n$ géné\-ri\-que\-ment inclus dans 
$X$ d'après le lemme \ref{l8}, qui est égal à un translaté d'un sous-groupe 
connexe contenant un borel d'après le théorème \ref{t1}.\ebew

\section{Commutateurs}
Dans un mauvais groupe de rang trois, Frécon \cite{Fr} trouve grâce aux
commutateurs un ensemble de Frécon de dimension deux ; il en conclut ensuite
que le groupe n'existe pas, par une méthode plus calculatoire que celle 
utilisée ici.
Dans la section précédente, nous avons généralisé la deuxième partie de sa 
démonstration au-delà de toute attente, mais pas la première, si bien que nous 
ne savons pas si notre Théorème \ref{t1} permet de montrer d'autres résultats 
d'inexistence. Dans cette section, nous montrerons que pour un groupe malsain 
$G$ avec borel $B$ abélien, on obtient au moins un $2\,RM(B)$-ensemble de 
Frécon, ce qui fournit une contradiction si $RM(G)=2\,RM(B)+1$, et en 
particulier pour $RM(G)=3$.

\satz\label{t2} Soit $G$ un groupe malsain avec borel $B$ abélien. Alors il y a 
un $2\,RM(B)$-ensemble de Frécon de rang strictement inférieur à $RM(G)$.\esatz
\bew Soit $g\not=1$ un commutateur non-trivial, et
$$X = \{ x\in G :\exists y\ [x,y] = g \}.$$
Si $x\in X$ et $y\in G$ est tel que $[x,y] = g$, alors $[x,y']=g$ pour tout 
$y'\in B(x)y$, et $xB(y')\subset X$.
Ces droites sont toutes différentes~: en effet, si $xB(y')=xB(y'')$ avec 
$y'\not= y''$ dans $B(x)$, alors $B(y')=B(y'')$ coupe $B(x)$ en deux points et 
lui est égal. Donc $g=[x,y']=1$, une contradiction.

Ainsi $X\subseteq\overline X^{\scriptscriptstyle 2\,RM(B)}$, et $RM(X)\ge 
2\,RM(B)$. Pour majorer le rang de $X$ on considère l'application 
$f:(x,y)\mapsto[x,y]$. L'image est invariant par conjugaison~; s'il ne contient 
qu'un nombre fini de classes de conjugaison, alors un commutateur $[x,y]$ de 
deux éléments génériques indépendants serait conjugué à son inverse $[y,x]$, ce 
qui est impossible. Par conséquent $RM(\im f)$ est strictement plus grand que 
le rang $RM(G)-RM(B)$ d'une classe de conjugaison, et un commutateur générique 
$g$ a une fibre $f^{-1}(g)$ de rang strictement inférieur à 
$2\,RM(G)-(RM(G)-RM(B))=RM(G)+RM(B)$. On en conclut que $RM(X)<RM(G)$.\ebew
Dans tout groupe malsain $G$ on a $RM(G)>2\,RM(B)$, car sinon pour $g\in G$ 
générique les deux double translatés $BgB$ et $Bg^{-1}B$ seraient génériques, 
donc d'intersection non-vide, donc égaux, ce qui fournirait une involution.
Si on fixe le rang du borel (abélien~!), le corollaire suivant exclut donc le 
cas minimal.
\kor Il n'y a pas de groupe $G$ malsain avec borel $B$ abélien de rang 
$RM(G)=2\,RM(B)+1$. En particulier, il n'y a pas de mauvais groupe de rang 
$3$.\ekor
\bew Dans un tel groupe, d'après le théorème \ref{t2} il y a un 
$2\,RM(B)$-ensemble de Frécon de rang $2\,RM(B)$~; d'après le corollaire 
\ref{k1} il y a un sous-groupe définissable connexe $H$ de rang $2\,RM(B)$ 
contenant un borel $B'$. Comme $H$ est malsain de borel $B'$, c'est 
absurde.\ebew

Par contre, on n'obtient pas de contradiction si $RM(G) = 3\,RM(B)$, même pour 
$RM(G)=6$ et $RM(B)=2$~: on aurait
$RM(X) = 5$, et les classes de conjugaison de commutateurs formeraient une
(toute petite !) famille de rang un. Quant à $RM(G)=4$ et $RM(B)=1$, ça 
voudrait dire aussi
que $RM(X) = 3$, mais que toutes les classes de conjugaison sauf un nombre
fini sont des commutateurs.


\begin{thebibliography}{9}

\bibitem{BP}A. V. Borovik et B. P. Poizat.
\newblock Simple groups of finite {M}orley rank without nonnilpotent
connected subgroups,
\newblock preprint deposited at VINITI, 1990.

\bibitem{Ch}Gregory Cherlin.
\newblock Groups of small Morley rank,
\newblock Ann.\ Math.\ Logic, 17(1--2):1--28, 1979.
 
\bibitem{Co}Luis Jaime Corredor.
\newblock Bad groups of finite Morley rank,
\newblock J. Symb.\ Logic, 54(3):768--773, 1989.

\bibitem{Fr}Olivier Fr\'econ.
\newblock Bad groups in the sense of Cherlin,
\newblock{\em preprint.}

\bibitem{Ne}Ali Nesin.
\newblock Nonsolvable groups of Morley rank $3$,
\newblock J. Algebra, 124(1):199--218, 1989.



\end{thebibliography}
\end{document}